
\input amstex
\documentstyle{amsppt}

\def\NS{NS_{\kappa\lambda}}
\def\Pkl{{\Cal P}_{\kappa}\lambda}
\def\cf{\operatorname{cf}}
\def\non{\operatorname{non}}
\def\cof{\operatorname{cof}}

\topmatter
\title
   Nowhere precipitousness of the non-stationary ideal over $\Pkl$
\endtitle
\author
   Yo Matsubara and Saharon Shelah
\endauthor
\thanks
The first author was partially supported by Grant-in-Aid for Scientific
Research (No.11640112), Ministry of Education, Science and Culture of Japan.
The second author was partially supported by The Israel Science Foundation
funded by the Israel Academy of Sciences and Humanities.
Publication [MsSh:758].
\endthanks

\abstract
   We prove that if $\lambda$ is a strong limit singular
   cardinal and $\kappa$ a regular uncountable cardinal $< \lambda$,
   then $\NS$, the non-stationary ideal over $\Pkl$, is nowhere
   precipitous.
   We also show that under the same hypothesis every stationary
   subset of $\Pkl$ can be partitioned into $\lambda^{< \kappa}$
   disjoint stationary sets.
\endabstract
\endtopmatter

\document
%
%
\head
   \S 1. Introduction
\endhead
Throughout this paper we let $\kappa$ denote an uncountable
regular cardinal and $\lambda$ a cardinal $\geq \kappa$.
Let $\NS$ denote the non-stationary ideal over $\Pkl$.
$\NS$ is the minimal $\kappa$-complete normal ideal
over $\Pkl$.
If $X$ is a stationary subset of $\Pkl$, then $\NS|X$
denotes the $\kappa$-complete normal ideal generated by
the members of $\NS$ and $\Pkl-X$.
We refer the reader to Kanamori [6, Section~25] for basic
facts about the combinatorics of $\Pkl$.
\par
Large cardinal properties of ideals have been investigated
by various authors.
One of the problems studied by these set theorists was
to determine which large cardinal properties can $\NS$
or $\NS|X$ bear for various $\kappa$, $\lambda$
and $X \subseteq \Pkl$.
In the course of this investigation, special interest
has been paid to two large cardinal properties, namely
precipitousness and saturation.
\par
If $\NS|X$ is not precipitous for every stationary $X \subseteq
\Pkl$, then we say that $\NS$ is {\it nowhere precipitous}.
In [8] Matsubara and Shioya proved that if $\lambda$ is a
strong limit singular cardinal and $\cf \lambda < \kappa$,
then $\NS$ is nowhere precipitous.
In \S 2 we extend this result by showing that $\NS$ is nowhere
precipitous if $\lambda$ is a strong limit singular cardinal.
\par
In [10] Menas conjectured the following:
\par
\proclaim{Menas' Conjecture}
   Every stationary subset of $\Pkl$ can be partitioned
   into $\lambda^{< \kappa}$ disjoint stationary sets.
\endproclaim
This conjecture implies that $\NS|X$ cannot
be $\lambda^{< \kappa}$-saturated for every
stationary $X\subseteq \Pkl$.
By the work of several set theorists we know that Menas'
Conjecture is independent of ZFC.
One of the most striking results concerning this conjecture
is the following theorem of Gitik [4].
\par
\proclaim{Gitik's Theorem}
   Suppose that $\kappa$ is a supercompact cardinal
   and $\lambda > \kappa$.
   Then there is a p.o. $\Bbb P$ that preserves
   cardinals $\geq \kappa$ such that
   $\Vdash_{\Bbb P}${\rm ``}$\exists X$ {\rm (}$X$ is
   a stationary subset of $\Pkl \wedge X$ cannot be
   partitioned into $\kappa^{+}$ disjoint stationary sets{\rm )''}.
\endproclaim
\par
In \S 2 we also show that if $\lambda$ is a strong limit
singular cardinal, then every stationary subset of $\Pkl$
can be partitioned into $\lambda^{< \kappa}$ disjoint
stationary sets.
Gitik [4] mentions that GCH fails in his model of a
``non-splittable'' stationary subset of $\Pkl$.
Our result shows that GCH {\it must} fail in such
a model of a non-splittable stationary subset of $\Pkl$
if $\lambda$ is singular.
\par
We often consider the poset ${\Bbb P}_{I}$ of $I$-positive
subsets of $\Pkl$ i.e. subsets of $\Pkl$ not belonging
to $I$, ordered by
$$
   X \leq_{{\Bbb P}_{I}} Y \Longleftrightarrow
   X \subseteq Y.
$$
We say that an ideal $I$ is ``proper'' if ${\Bbb P}_{I}$
is a proper poset.
In [9] Matsubara proved the following result:
\par
\proclaim{Proposition}
   Let $\delta$ be a cardinal $\geq 2^{2^{2^{\lambda}}}$.
   If there is a ``proper'' $\lambda^{+}$-complete normal
   ideal over ${\Cal P}_{\lambda^{+}}\delta$
   then $NS_{\aleph_{1}\lambda}$ is precipitous.
\endproclaim
\par
It is not known whether $\NS$ can be precipitous for
singular $\lambda$.
In [1] it is conjectured that $\NS$ cannot be precipitous
if $\lambda$ is singular.
Therefore it is interesting to ask the following question:
\par
\proclaim{Question}
   Can $\Pkl$ bear a ``proper'' $\kappa$-complete normal
   ideal where $\kappa$ is the successor cardinal of
   a singular cardinal?
\endproclaim
\par
In \S 3 we give a negative answer to this question.
\par
%
%
\head
   \S 2. On $\NS$ for strong limit singular $\lambda$
\endhead
We first state our main results.
\par
\proclaim{Theorem 1}
   If $\lambda$ is a strong limit singular cardinal,
   then $\NS$ is nowhere precipitous.
\endproclaim
\par
\proclaim{Theorem 2}
   If $\lambda$ is a strong limit singular cardinal,
   then every stationary subset of $\Pkl$ can be
   partitioned into $\lambda^{< \kappa}$ disjoint
   stationary sets.
\endproclaim
\par
One of the key ingredients of our proof of the main
results is Lemma 3.
  Part (ii) of Lemma 3 was proved in Matsubara [7].
  Part (i) appeared in Matsubara-Shioya [8].
For the proof of Part (ii) we refer the reader to
Kanamori [6, page 345].
However we will present the proof of (i) because
the idea of this proof will be used later.
\par
\proclaim{Lemma 3}
   If $2^{< \kappa} < \lambda^{< \kappa} = 2^{\lambda}$,
   then
   \par
   {\rm (i)} $\NS$ is nowhere precipitous
   \par
   {\rm (ii)} every stationary subset of $\Pkl$ can be
   partitioned into $\lambda^{< \kappa}$ disjoint
   stationary sets.
\endproclaim
\par
Before we present the proof of part (i), we make
some comments concerning this lemma.
First note that the hypothesis of our lemma is
satisfied if $\lambda$ is a strong limit cardinal
with $\cf\lambda < \kappa$.
Secondly under this hypothesis every unbounded subset
of $\Pkl$ must have a size of $2^{\lambda}$.
We also note that Lemma~3 can be generalized in the
following manner:
\par
For an ideal $I$ over some set $A$, we let $\non(I) =
\min \{ |X| \mid X \subseteq A,\,X \notin I \}$
and $\cof(I) = \min \{ |J| \mid J \subseteq I,\,
\forall X \in I,\,\exists Y \in J\,(X \in Y)\}$.
The proof of Lemma~3 actually shows that
if $\non(I) = \cof(I)$ then $I$ is nowhere precipitous
(i.e. for every $I$-positive $X$, $I|X$ is not precipitous)
and every $I$-positive subset $X$ of $A$ can be
partitioned into $\non(I)$ many disjoint $I$-positive
sets.
\par
\demo{Proof of Lemma~{\rm 3~(i)}}
   For $I$ an ideal over $\Pkl$, let $G(I)$ denote
   the following game between two players, Nonempty
   and Empty:
   Nonempty and Empty alternately choose $I$-positive
   sets $X_{n}$, $Y_{n} \subseteq \Pkl$ respectively
   so that $X_{n} \supseteq Y_{n} \supseteq X_{n+1}$
   for $n = 1, 2, \dots$.
   After $\omega$ moves, Empty wins $G(I)$
   if $\bigcap_{n \in \omega-\{0\}} X_{n} = \emptyset$.
   See [3] for a proof of the following characterization.
   \par
   \proclaim{Proposition}
     $I$ is nowhere precipitous if and only if Empty
     has a winning strategy in $G(I)$.
   \endproclaim
   \par
   Let $\langle f_{\alpha} \mid \alpha < 2^{\lambda} \rangle$
   enumerate functions from $\lambda^{< \omega}$
   into $\Pkl$.
   For a function $f : \lambda^{< \omega} \to \Pkl$,
   we let $C(f) = \{ s \in \Pkl \mid
   \bigcup f'' s^{< \omega} \subseteq s \}$.
   For $X \subseteq \Pkl$, $X$ is stationary if and
   only if $C(f_{\alpha}) \cap X \ne \emptyset$ for
   every $\alpha < 2^{\lambda}$.
   \par
   We now describe Empty's strategy in $G(\NS)$ using
   the hypothesis $2^{< \kappa} < \lambda^{< \kappa} = 2^{\lambda}$.
   Suppose that $X_{1}$ is Nonempty's first move.
   Choose $\langle s^{1}_{\alpha} \mid \alpha < 2^{\lambda} \rangle$,
   a sequence of elements of $X_{1}$ by induction on $\alpha$
   in the following manner:
   Let $s^{1}_{0}$ be any element of $X_{1} \cap C(f_{0})$.
   Suppose we have $\langle s^{1}_{\alpha} \mid \alpha < \beta \rangle$
   for some $\beta < 2^{\lambda}$.
   Since $\{ s^{1}_{\alpha} \mid \alpha < \beta \}$ is
   a non-stationary, in fact bounded, subset of $\Pkl$,
   $X_{0} - \{ s^{1}_{\alpha} \mid \alpha < \beta \}$ is
   stationary.
   Pick an element from $(X_{0} - \{ s^{1}_{\alpha} \mid
   \alpha < \beta \}) \cap C(f_{\beta})$ and call it $s^{1}_{\beta}$.
   Let Empty play $Y_{1} = \{ s^{1}_{\alpha} \mid \alpha < 2^{\lambda} \}$.
   It is easy to see that $Y_{1}$ is a stationary subset
   of $\Pkl$.
   Inductively suppose Nonempty plays his $n+1$-st
   move $X_{n+1}$ immediately following Empty's $n$-th
   move $Y_{n} = \{ s_{\alpha}^{n} \mid \alpha < 2^{\lambda} \}$.
   Choose $\langle s_{\alpha}^{n+1} \mid \alpha < 2^{\lambda} \rangle$,
   a sequence from $X_{n+1}$ in the following manner:
   Let $s_{0}^{n+1}$ be any element of $(X_{n+1}-\{ s_{0}^{n} \})
   \cap C(f_{0})$.
   Suppose we have $\langle s_{\alpha}^{n+1} \mid \alpha < \beta \rangle$,
   for some $\beta < 2^{\lambda}$.
   Pick an element of the stationary set $(X_{n+1} \cap C(f_{\beta})) -
   (\{ s_{\alpha}^{n+1} \mid \alpha < \beta \} \cap
   \{ s_{\alpha}^{n} \mid \alpha \leq \beta \})$ and
   call it $s_{\beta}^{n+1}$.
   Let Empty play $Y_{n+1} = \{ s_{\alpha}^{n+1} \mid
   \alpha < 2^{\lambda} \}$.
   This defines a strategy for Empty.
   \par
   \proclaim{Claim}
     The strategy described above is a winning strategy for Empty.
   \endproclaim
   \par
   \demo{Proof of Claim}
     Suppose $X_{1}, Y_{1}, X_{2}, Y_{2}, \dots$ is a run of
     the game $G(\NS)$ where Empty followed the above strategy.
     We want to show that $\bigcap_{n \in \omega-\{0\}} Y_{n} =
     \emptyset$.
     Suppose otherwise.
     Let $t$ be an element of $\bigcap_{n \in \omega-\{0\}} Y_{n}$.
     Then for each $m \in \omega - \{ 0 \}$, there is a unique
     ordinal $\alpha_{m} < 2^{\lambda}$ such that $s_{\alpha_{m}}^{m} = t$.
     But by the way the $s_{\alpha}^{n}$s are chosen, $s_{\alpha_{0}}^{0} =
     s_{\alpha_{1}}^{1} = s_{\alpha_{2}}^{2} = \cdots$
     implies $\alpha_{0} > \alpha_{1} > \alpha_{2} > \cdots$.
     This is impossible.
     Thus we must have $\bigcap_{n \in \omega-\{0\}} Y_{n} =
     \emptyset$.
     \qed
   \enddemo
   \hfill End of proof of Lemma~3~(i). \qed
\enddemo
\par
We now prove Theorem~2 using Lemma~3 and Theorem~1.
\par
\demo{Proof of Theorem~{\rm 2}}
   Let $\lambda$ be a strong limit singular cardinal.
   If $\cf\lambda < \kappa$ then by Lemma~3~(ii), we are done.
   So assume $\cf\lambda \geq \kappa$.
   In this case we have $\lambda^{< \kappa} = \lambda$.
   Therefore it is enough to show that $\NS|X$ is
   not $\lambda$-saturated for every stationary $X \subseteq \Pkl$.
   But this is a consequence of $\NS$ being nowhere
   precipitous.
   In fact we know that $\NS|X$ cannot be $\lambda^{+}$-saturated
   for every stationary $X \subseteq \Pkl$.
   \qed
\enddemo
\par
We need some preparation to present the proof of Theorem~1.
Let $\lambda$ be a strong limit singular cardinal
and $\kappa$ be a regular uncountable cardinal $< \lambda$.
If $\cf\lambda < \kappa$ then by Lemma~3 we conclude
that $\NS$ is nowhere precipitous.
\par
 From now on let us assume that $\lambda$ is a strong limit
cardinal with $\kappa \leq \cf\lambda < \lambda$.
Let $\langle \lambda_{\alpha} \mid \alpha < \cf\lambda \rangle$
be a continuous increasing sequence of strong limit
singular cardinals converging to $\lambda$.
The following lemma is another key ingredient of our proof.
\par
\proclaim{Lemma~4}
   For every $X \subseteq \Pkl$, if for each $\alpha < \cf\lambda$
   with $\cf\alpha < \kappa$, $|\{ t \in X \mid
   \sup(t) = \lambda_{\alpha} \}| < 2^{\lambda_{\alpha}}$,
   then $X$ is non-stationary.
\endproclaim
\par
\demo{Proof of Lemma~{\rm 4}}
   Since $\{ t \in X \mid \sup(t) \notin t \}$ is a club subset
   of $\Pkl$, without loss of generality we may assume
   that $\sup(t) \notin t$ for every $t$ in $X$.
   For each $\alpha < \cf\lambda$ with $\cf\alpha < \kappa$,
   we let $X_{\alpha} = \{ t \in X \mid \sup(t) = \lambda_{\alpha} \}$.
   We need the following fact from pcf theory by S.~Shelah.
   \par
   \proclaim{Fact}
     There is a club subset $C \subseteq \cf\lambda$
     such that $\operatorname{pp}(\lambda_{\alpha}) =
     2^{\lambda_{\alpha}}$ for every $\alpha \in C$.
   \endproclaim
   \par
   We refer the reader to Shelah [11] and Holz-Steffens-Weitz
   [5] for the pcf theory needed.
   The proof of the above fact can be obtained by combining
   Conclusion 5.13 [11, page 414] and Theorem 9.1.3 [5, page 271].
   \par
   For each $\alpha \in C$ with $\cf\alpha < \kappa$,
   let $a_{\alpha}$ be a set of regular cardinals cofinal
   in $\lambda_{\alpha}$ such that
   \par
   \hskip1em (a) every member of $a_{\alpha}$ is above $\cf\lambda$
   \par
   \hskip1em (b) $|a_{\alpha}| = \cf\lambda_{\alpha}$, and
   \par
   \hskip1em (c) $\exists \delta_{\alpha} > |X_{\alpha}|$
   $[\delta_{\alpha} \in \operatorname{pcf}(a_{\alpha})]$
   \par
   Let $a = \bigcup \{ a_{\alpha} \mid \alpha \in C \wedge
   \cf\alpha < \kappa \}$.
   Let $\langle f_{\beta} \mid \beta < \lambda \rangle$
   enumerate all of the members of $\{ f \mid$ $f$ is
   a function, $\operatorname{domain}(f)$ is a bounded subset
   of $\lambda$, and $f$ is regressive i.e. $f(\gamma) < \gamma$
   for every $\gamma \in \operatorname{domain}(f) \}$.
   \par
   For each $t \in \Pkl$ we define $g_{t} \in \prod a$
   by letting $g_{t}(\sigma) = \sup \{ f_{\beta}(\sigma)+1 \mid
   \beta \in t \wedge \sigma \in \operatorname{dom}(f_{\beta}) \}$,
   if $\sigma \in \bigcup_{\beta \in t} \operatorname{domain}(f_{\beta})$,
   and $g_{t}(\sigma) = 0$ otherwise.
   Note that $|t| < \kappa \leq \cf\lambda < \min(a)$
   guarantees $g_{t} \in \prod a$.
   Now by (c) in the definition of $a_{\alpha}$s and the fact
   that $\{ g_{t}\upharpoonright a_{\alpha} \mid t \in X_{\alpha} \}$ is a 
subset
   of $\prod a_{\alpha}$ of cardinality $\leq |X_{\alpha}| <
   \delta_{\alpha} \in \operatorname{pcf}(a_{\alpha})$,
   there is some $h_{\alpha} \in \prod a_{\alpha}$ such
   that $\forall t \in X_{\alpha}$ $[g_{t}\upharpoonright a_{\alpha}
   <_{J_{<\delta_{\alpha}}(a_{\alpha})} h_{\alpha}]$.
   Therefore
   $$
     \forall t \in X_{\alpha}\
     \exists \sigma \in a_{\alpha}\
     [g_{t}(\sigma) < h_{\alpha}(\sigma)]
     \tag{1}
   $$
   holds.
   As $\min(a) > \cf\lambda$ and $a = \bigcup \{ a_{\alpha} \mid
   \alpha \in C \wedge \cf\alpha < \kappa \}$,
   there is $h \in \prod a$ such that $h_{\alpha} <
   h\upharpoonright a_{\alpha}$
   for every $\alpha \in C$ with $\cf\alpha < \kappa$.
   \par
   Let $W = \{ t \in \Pkl \mid$ (i) for some $\alpha \in C$
   $\sup(t) = \lambda_{\alpha}$ with $\cf\alpha < \kappa$, and
   (ii) if $\delta \in t$ then for some $\beta \in t$,
   $h\upharpoonright (a \cap \delta) = f_{\beta} \}$.
   Note that $W$ is a club subset of $\Pkl$.
   \par
   \proclaim{Claim}
     $X \cap W = \emptyset$.
   \endproclaim
   \par
   \demo{Proof of Claim}
     Suppose otherwise, say $t \in X \cap W$.
     By (i) in the definition of $W$, $t \in X_{\alpha}$
     for some $\alpha \in C$ with $\cf\alpha < \kappa$.
     By (1) we have
     $$
       \exists \sigma \in a_{\alpha}\
       [g_{t}(\sigma) < h_{\alpha}(\sigma)].
       \tag{2}
     $$
     Since $\sup(t) = \lambda_{\alpha}$, there must be
     some $\delta \in t$ such that $\delta > \sigma$.
     Now by (ii) in the definition of $W$, $h\upharpoonright
     (a \cap \delta) =
     f_{\beta}$ for some $\beta \in t$.
     Since $\sigma \in a \cap \delta$, $h(\sigma) = f_{\beta}(\sigma)$.
     By the definition of $g_{t}$ we have $f_{\beta}(\sigma) <
     g_{t}(\sigma)$.
     From $h_{\alpha} < h\upharpoonright a_{\alpha}$, we know
     $h_{\alpha}(\sigma) <
     h(\sigma)$.
     Therefore we have $h_{\alpha}(\sigma) < g_{t}(\sigma)$
     contradicting (2).
     \qed
   \enddemo
   \hfill End of proof of Lemma~4. \qed
\enddemo
\par
For each $\alpha < \cf\lambda$ with $\cf \alpha < \kappa$,
let us fix a sequence $\langle f_{\xi}^{\alpha} \mid
\xi < 2^{\lambda_{\alpha}} \rangle$ that enumerates members
of $\{ f \mid f$ is a function such that $\operatorname{domain}(f) \subseteq
\lambda_{\alpha}^{< \omega}$ and $\operatorname{range}(f) \subseteq
\lambda_{\alpha} \}$.
Furthermore for each function $f$
with $\operatorname{domain}(f) \subseteq \lambda_{\alpha}^{< \omega}$
and $\operatorname{range}(f) \subseteq \lambda_{\alpha}$,
we let $C_{\alpha}[f] = \{ t \in \Pkl \mid t^{< \omega} \subseteq
\operatorname{domain}(f)$, $\sup(t) = \lambda_{\alpha}$, and $t$ is
closed under $f \}$.
We need the following lemma to present the proof of Theorem~1.
\par
\proclaim{Lemma~5}
   Suppose $X$ is a stationary subset of $\Pkl$.
   For every $Y \subseteq \{ s \in \Pkl \mid s \cap \kappa \in \kappa \}$,
   if for each $\alpha < \cf\lambda$ with $\cf\alpha < \kappa$
   the following condition $(*)$ holds, then $Y$ is stationary.
   $$
     \forall \xi < 2^{\lambda_{\alpha}}\
     (|C_{\alpha}[f_{\xi}^{\alpha}] \cap X| = 2^{\lambda_{\alpha}}
     \longrightarrow
     C_{\alpha}[f_{\xi}^{\alpha}] \cap Y \ne \emptyset)
     \tag{$*$}
   $$
\endproclaim
\par
\demo{Proof of Lemma~{\rm 5}}
   Since $s \cap \kappa \in \kappa$ for every $s \in Y$, to show
   that $Y$ is stationary it is enough to show that $Y \cap C[g] \ne
   \emptyset$ for every function $g:\lambda^{< \omega} \to \lambda$
   where $C[g]$ denotes the set $\{ t \in \Pkl \mid
   g''t^{< \omega} \subseteq t \}$.
   For the proof of this fact, we refer the reader to
   Foreman-Magidor-Shelah [2, Lemma~0].
   Let us fix a function $g:\lambda^{< \omega} \to \lambda$.
   Now we let $E = \{ \alpha < \cf\lambda \mid \cf\alpha < \kappa \}$
   and for each $\alpha \in E$ we let $W_{\alpha} = \{ s \in \Pkl \mid
   \sup(s) = \lambda_{\alpha} \wedge \lambda_{\alpha} \notin s \}$.
   Note that $\bigcup_{\alpha \in E} W_{\alpha}$ is a club subset
   of $\Pkl$.
   For each $\alpha \in E$, we let $g_{\alpha}$
   denote $g \cap (\lambda_{\alpha}^{< \omega} \times \lambda_{\alpha})$.
   Now partition $E$ into two sets $E^{+}$ and $E^{-}$ where
   $$
     \align
       E^{+} &= \{ \alpha \in E \mid |C_{\alpha}[g_{\alpha}] \cap X| =
         2^{\lambda_{\alpha}} \} \quad \text{and} \cr
       E^{-} &= \{ \alpha \in E \mid |C_{\alpha}[g_{\alpha}] \cap X| <
         2^{\lambda_{\alpha}} \}.
     \endalign
   $$
   \par
   We need the following:
   \par
   \proclaim{Claim}
     $X \cap \bigcup \{ W_{\alpha} \mid \alpha \in E^{-} \}$ is
     non-stationary.
   \endproclaim
   \par
   \demo{Proof}
     It is enough to show that $Z = C[g] \cap X \cap
     \bigcup \{ W_{\alpha} \mid \alpha \in E^{-} \}$ is
     non-stationary.
     Note that for each $\alpha \in E^{+}$,
     $Z \cap W_{\alpha} = \emptyset$ and for each $\alpha \in E^{-}$,
     $Z \cap W_{\alpha} \subseteq C_{\alpha}[g_{\alpha}] \cap X$.
     Therefore $|Z \cap W_{\alpha}| < 2^{\lambda_{\alpha}}$ for
     every $\alpha \in E$.
     Hence, by Lemma~4, we conclude that $Z$ is non-stationary.
     \qed
   \enddemo
   \par
   From Claim we know that $X \cap \bigcup
   \{ W_{\alpha} \mid \alpha \in E^{+} \}$ is stationary.
   Pick an element $\alpha^{*}$ from $E^{+}$.
   Consider the partial function $g_{\alpha^{*}}$ $(= g \cap
   (\lambda_{\alpha^{*}}^{< \omega} \times \lambda_{\alpha^{*}}))$.
   Let $\xi^{*} < 2^{\lambda_{\alpha^{*}}}$ be
   such that $f_{\xi^{*}}^{\alpha^{*}} = g_{\alpha^{*}}$.
   Since $\alpha^{*} \in E^{+}$,
   we have $|C_{\alpha^{*}}[g_{\alpha^{*}}] \cap X| =
   2^{\lambda_{\alpha^{*}}}$.
   Since $f_{\xi^{*}}^{\alpha^{*}} = g_{\alpha^{*}}$
   and $Y$ satisfies condition $(*)$, we know
   that $C_{\alpha^{*}}[g_{\alpha^{*}}] \cap Y \ne \emptyset$.
   Therefore $C[g] \cap Y \ne \emptyset$ showing that $Y$ is
   stationary.
   \par
   \hfill End of proof of Lemma~5. \qed
\enddemo
\par
Finally we are ready to complete the proof of Theorem~1.
To present a winning strategy for Empty in the game $G(\NS)$,
we introduce some new types of games.
For each $\alpha \in E = \{ \alpha < \cf\lambda \mid
\cf\alpha < \kappa \}$, we define the game $G_{\alpha}$ between
Nonempty and Empty as follows:
Nonempty and Empty alternately choose sets $X_{n}, Y_{n} \subseteq
W_{\alpha} = \{ s \in \Pkl \mid \sup(s) = \lambda_{\alpha} \notin s \}$
respectively so that $X_{n} \supseteq Y_{n} \supseteq X_{n+1}$
and $\forall \xi < 2^{\lambda_{\alpha}}$ $(|C_{\alpha}[f_{\xi}^{\alpha}] \cap
X_{n}| = 2^{\lambda_{\alpha}} \longrightarrow C_{\alpha}[f_{\xi}^{\alpha}]
\cap Y_{n} \ne \emptyset)$ for $n = 1, 2, \dots$.
Empty wins $G_{\alpha}$ iff $\bigcap_{n \in \omega-\{0\}} Y_{n} =
\emptyset$.
\par
By the same argument as the proof of Lemma~3~(i), we know
that Empty has a winning strategy, say $\tau_{\alpha}$,
in the game $G_{\alpha}$ for each $\alpha \in E$.
Now we show how to combine the strategies $\tau_{\alpha}$s
to produce a winning strategy for Empty in $G(\NS)$.
Suppose $X_{1}$ is Nonempty's first move in $G(\NS)$.
We let $X_{1}^{*} = X_{1} \cap \{ s \in \Pkl \mid
s \cap \kappa \in \kappa \} \cap
\bigcup \{ W_{\alpha} \mid \alpha \in E \}$.
Since $\{ s \in \Pkl \mid s \cap \kappa \in \kappa \} \cap
\bigcup \{ W_{\alpha} \mid \alpha \in E \}$ is a club subset
of $\Pkl$, $X_{1}^{*}$ is stationary in $\Pkl$.
For each $\alpha \in E$, we simulate a run of the game $G_{\alpha}$
as follows:
Let us pretend that Nonempty's first move in $G_{\alpha}$
is $X_{1}^{*} \cap W_{\alpha}$.
Let Empty play her strategy $\tau_{\alpha}$, so Empty's first
move is $\tau_{\alpha}(\langle X_{1}^{*} \cap W_{\alpha} \rangle)$.
Now in the game $G(\NS)$, let Empty play $Y_{1} =
\bigcup \{ \tau_{\alpha}(\langle X_{1}^{*} \cap W_{\alpha} \rangle) \mid
\alpha \in E \}$.
Lemma~5 guarantees that $Y_{1}$ is stationary in $\Pkl$.
In general if $\langle X_{1}^{*}, Y_{1}, X_{2}, Y_{2}, \dots, X_{n} \rangle$
is a run of $G(\NS)$ up to Nonempty's $n$-th move,
then we let Empty play $Y_{n} = \bigcup \{ \tau_{\alpha}
(\langle X_{1}^{*} \cap W_{\alpha}, X_{2} \cap W_{\alpha}, \dots,
X_{n} \cap W_{\alpha} \rangle) \mid \alpha \in E \}$.
Once again we know $Y_{n}$ is a stationary subset of $X_{n}$.
For each $\alpha \in E$, since $\tau_{\alpha}$ is a winning strategy
in
$G_{\alpha}$  we have
$$
   \bigcap_{n \in \omega-\{0\}}
   \tau_{\alpha}(\langle X_{1}^{*} \cap W_{\alpha},
   X_{2} \cap W_{\alpha}, \dots, X_{n} \cap W_{\alpha} \rangle) =
   \emptyset.
$$
Because the $W_{\alpha}$s are pairwise disjoint,
we conclude that $\bigcap_{n \in \omega-\{0\}} Y_{n} = \emptyset$.
Therefore we have a winning strategy for Empty in
the game $G(\NS)$.
This proves that $\NS$ is nowhere precipitous for every
strong limit singular $\lambda$.
\par
\hfill End of proof of Theorem~1. \qed
\par
%
%
\head
  \S 3.  On ``proper'' ideals over $\Pkl$
\endhead
First we define that we mean by a ``proper'' ideal.
\par
\demo{Definition}
   An ideal $I$ over a set $A$ is a ``proper'' ideal if
   the corresponding p.o. ${\Bbb P}_{I}$ is proper (in the sense
   of proper forcing).
\enddemo
\par
We refer the reader to Shelah [12] for the background of properness.
\par
As we mentioned in \S 1, we are interested in the question of
whether it is possible to have a $\kappa$-complete normal ``proper''
ideal over $\Pkl$ where $\kappa$ is the successor of some singular
cardinal.
We give a negative answer to this question.
Here we present a more general result.
\par
\proclaim{Theorem~6}
   {\rm (i)} Suppose $I$ is a $\kappa$-complete normal ideal
   over $\kappa$.
   If $\{ \alpha < \kappa \mid \cf\alpha = \delta \} \notin I$
   for some cardinal $\delta$ satisfying $\delta^{+} < \kappa$,
   then $I$ is not ``proper''.
   \par
   {\rm (ii)} Suppose $I$ is a $\kappa$-complete normal ideal
   over $\Pkl$.
   If $\{ s \in \Pkl \mid \cf(s \cap \kappa) = \delta \} \notin I$
   for some cardinal $\delta$ satisfying $\delta^{+} < \kappa$,
   then $I$ is not ``proper''.
\endproclaim
\par
Note that if $\kappa$ is the successor cardinal of a singular
cardinal, then every $\kappa$-complete normal ideal
over $\Pkl$ satisfies the hypothesis of (ii).
\par
\demo{Proof of Theorem~{\rm 6}}
   Since the proof of (ii) is identical to that of (i),
   we only present the proof of (i).
   \par
   Let $I$ and $\delta$ be as in the hypothesis of (i).
   First note that if $\delta = \aleph_{0}$ then
   the set $\{ \alpha < \kappa \mid \cf\alpha = \delta \}$
   forces ``$\cf\kappa = \aleph_{0}$'' showing ${\Bbb P}_{I}$
   cannot be proper.
   Therefore we may assume that $\delta$ is uncountable.
   \par
   We need the following claim:
   \par
   \proclaim{Claim~1}
     There are a stationary subset $E$ of $\{ \alpha < \kappa \mid
     \cf\alpha = \aleph_{0} \}$ and an $I$-positive subset $X$
     of $\{ \alpha < \kappa \mid \cf\alpha = \delta \}$
     such that $E \cap \alpha$ is non-stationary for every $\alpha$
     in $X$.
   \endproclaim
   \par
   \demo{Proof}
     Let $\{ E_{\gamma} \mid \gamma < \delta^{+} \}$ be a family
     of pairwise disjoint stationary subsets of $\{ \alpha < \kappa \mid
     \cf\alpha = \aleph_{0} \}$.
     For each $\alpha < \kappa$ with $\cf\alpha = \delta$,
     there must be a club subset of $\alpha$ with cardinality $\delta$.
     Therefore for such an ordinal $\alpha$, there is
     some $\gamma_{\alpha} < \delta^{+}$ such
     that $E_{\gamma_{\alpha}} \cap \alpha$ is non-stationary.
     By the $\kappa$-completeness of $I$,
     there is some $\gamma^{*} < \delta^{+}$
     such that $X = \{ \alpha < \kappa \mid \cf\alpha = \delta \wedge
     \gamma_{\alpha} = \gamma^{*} \}\notin I$.
     If we let $E = E_{\gamma^{*}}$, then $E \cap \alpha$ is
     non-stationary for every $\alpha$ in $X$.
     \qed
   \enddemo
   \par
   For each $\alpha$ from $X$, let $c_{\alpha}$ be a club subset
   of $\alpha$ with $c_{\alpha} \cap E = \emptyset$.
   Let $\vec{C}$ denote $\langle c_{\alpha} \mid \alpha \in X \rangle$.
   Let $\chi$ be a large enough regular cardinal.
   Assume that $N$ is a countable elementary substructure
   of $\langle H(\chi), \epsilon \rangle$
   satisfying $\{ I, E, X, \vec{C} \} \subseteq N$
   and $\sup(N \cap \kappa) \in E$.
   \par
   We are ready to show that $I$ is not ``proper''.
   \par
   \proclaim{Claim~2}
     If $Y$ is a subset of $X$ such that $Y \notin I$
     {\rm (}therefore $Y \in {\Bbb P}_{I}$ and $Y \leq X${\rm )},
     then $Y$ is not $(N, {\Bbb P}_{I})$-generic.
   \endproclaim
   \par
   Claim~2 implies that ${\Bbb P}_{I}$ is not proper.
   \par
   \demo{Proof of Claim~{\rm 2}}
     Suppose otherwise.
     Assume that there exists $Y \leq X$ such that $Y$
     is $(N, {\Bbb P}_{I})$-generic.
     \par
     For each $\alpha < \kappa$ we define a function $f_{\alpha}:
     X \to \kappa$ by $f_{\alpha}(\gamma) =
     \operatorname{Min}(c_{\gamma}-\alpha)$
     if $\gamma > \alpha$, and $f_{\alpha}(\gamma) = 0$ otherwise.
     It is clear that $f_{\alpha} \in N$ for each $\alpha \in N \cap \kappa$.
     \par
     For each $\alpha \leq \beta < \kappa$,
     we let $T_{\beta}^{\alpha} = \{ \gamma \in X \mid
     f_{\alpha}(\gamma) = \beta \}$.
     For each fixed $\alpha < \kappa$, using the normality of $I$,
     we see that $\{ T_{\beta}^{\alpha} \mid
     \alpha \leq \beta < \kappa, T_{\beta}^{\alpha} \notin I \}$ is
     a maximal antichain below $X$ in ${\Bbb P}_{I}$.
     Let $\vec{T}^{\alpha} = \langle T_{\beta}^{\alpha} \mid
     \alpha \leq \beta <\kappa, T_{\beta}^{\alpha} \notin I \rangle$.
     It is clear that $\vec{T}^{\alpha} \in N$
     for $\alpha \in N \cap \kappa$.
     \par
     Since $Y$ is $(N, {\Bbb P}_{I})$-generic,
     for $\alpha \in N \cap \kappa$ $\{ T_{\beta}^{\alpha} \mid
     \alpha \leq \beta \wedge \beta \in N \cap \kappa \wedge
     T_{\beta}^{\alpha} \notin I \}$ is predense below $Y$
     in ${\Bbb P}_{I}$.
     So we must have $Y - \bigcup \{ T_{\beta}^{\alpha} \mid
     \alpha \leq \beta \wedge \beta \in N \cap \kappa \wedge
     T_{\beta}^{\alpha} \notin I \} \in I$ for
     each $\alpha \in N \cap \kappa$.
     Let $Y_{\alpha} = Y - \bigcup \{ T_{\beta}^{\alpha} \mid
     \alpha \leq \beta \wedge \beta \in N \cap \kappa \wedge
     T_{\beta}^{\alpha} \notin I \}$.
     We have $\bigcup_{\alpha \in N \cap \kappa} Y_{\alpha} \in I$.
     This implies $Y - \bigcup_{\alpha \in N \cap \kappa}
     Y_{\alpha} \notin I$.
     Let $\gamma^{*}$ be an element
     of $Y - \bigcup_{\alpha \in N \cap \kappa} Y_{\alpha}$
     with $\gamma^{*} > \sup(N \cap \kappa)$.
     Note that $\gamma^{*} \in
     Y - Y_{\alpha}$ for each $\alpha \in N \cap \kappa$.
     Hence if $\alpha \in N \cap \kappa$,
     then there exists $\beta_{\alpha} \in N \cap \kappa$
     such that $\gamma^{*} \in T_{\beta}^{\alpha}$.
     Thus $f_{\alpha}(\gamma^{*}) = \beta_{\alpha} \in N \cap \kappa$
     for each $\alpha \in N \cap \kappa$.
     This means that $\operatorname{Min}(c_{\gamma^{*}}-\alpha) \in
     N \cap \kappa$ for each $\alpha \in N \cap \kappa$,
     showing $c_{\gamma^{*}} \cap N$ is unbounded
     in $\sup(N \cap \kappa)$.
     \par
     Since $\sup(N \cap \kappa) < \gamma^{*}$,
     we must have $\sup(N \cap \kappa) \in c_{\gamma^{*}}$.
     But this implies $\sup(N \cap \kappa) \in c_{\gamma^{*}} \cap E$
     which contradicts $c_{\alpha} \cap E = \emptyset$
     for each $\alpha \in X$ and $\gamma^{*} \in Y \subseteq X$.
     This contradiction shows that $Y$ cannot
     be $(N, {\Bbb P}_{I})$-generic.
     \qed
   \enddemo
   \par
   \hfill End of proof of Theorem~6. \qed
\enddemo
\par
%
%

\enddocument